\documentclass{amsart}
\nonstopmode
\usepackage{graphics}
\usepackage{epsfig}
\usepackage[frame,matrix,curve,arrow,dvips]{xy}
\usepackage{amssymb}
\usepackage{amsmath}
\usepackage{amsthm}

\newtheorem{thm}{Theorem}[section]
\newtheorem{lem}[thm]{Lemma}
\newtheorem{pro}[thm]{Proposition}

\newtheorem{cor}[thm]{Corollary}
\newtheorem{rem}[thm]{Remark}

\def\hpic #1 #2 {\mbox{$\begin{array}[c]{l} \epsfig{file=#1,height=#2} \end{array}$}}
\def\vpic #1 #2 {\mbox{$\begin{array}[c]{l} \epsfig{file=#1,width=#2} \end{array}$}}

\def\E #1 #2{{ _{#1} E _{#2} \ }}
\def\A #1 #2{{ A _{#1} ^{#2} \ }}
\def\e #1 #2{{ _{#1} e _{#2} }}

\def\ss{\subset}

\def\E{\mbox{$\exists$}}
\def\A{\mbox{$\forall$}}

\DeclareMathOperator{\aut}{Aut}

\newcommand{\dgre}{\lambda}
\newcommand{\dgge}{\rho}
\newcommand{\gre}{\nu}
\newcommand{\gge}{\theta}
\newcommand{\myA}{\Lambda} 
\newcommand{\myB}{{\rm P}}
\newcommand{\myC}{{\rm N}}
\newcommand{\myD}{\Theta}
\newcommand{\grp}{\Gamma}

\newcommand{\dgrp}{\grp'}
\newcommand{\alg}{A}
\newcommand{\salg}{B}
\newcommand{\ialg}{C}
\newcommand{\palg}{C}
\newcommand{\calg}{D}
\newcommand{\grph}{\mathcal{G}}

\begin{document}
\title{A note on intermediate subfactors of Krishnan-Sunder subfactors}
\author{Bina Bhattacharyya} 
\address{ Deephaven \\
          4699 Old Ironsides Dr., \#210\\
	  Santa Clara, CA 94054 \\
	 US}
\email{Bina\_Bhattacharyya\_91@post.harvard.edu}
\begin{abstract}
A Krishnan-Sunder subfactor $R_U \ss R$ of index $k^2$ is constructed
from a permutation biunitary matrix $U\in M_p(\mathbb{C})\otimes
M_k(\mathbb{C})$, i.e. the entries of $U$ are either 0 or 1 and both
$U$ and its block transpose are unitary.  The author previously showed
that every irreducible Krishnan-Sunder subfactor has an intermediate
subfactor by exhibiting the associated Bisch projection.  The author
has also shown in a separate paper that the principal and dual graphs
of the intermediate subfactor are the same as those of the subfactor
$R^{\grp} \ss R^{H}$, where $H\ss \grp$ is an inclusion of finite groups with an
outer action on $R$.  In this paper we give a direct proof that the
intermediate subfactor is isomorphic to $R^{\grp} \ss R^{H}$.
\end{abstract}
\maketitle

\section{Background and Introduction}
\label{sec:krishn-sund-subf}
There is a well-known way of constructing subfactors of the
hyperfinite ${\rm II}_1$ subfactor $R$ from certain squares of
finite-dimensional $C^\ast$--algebras algebras.  Since we will use
this construction repeatedly, we review it briefly.  Suppose we have a
square of finite-dimensional $C^\ast$--algebras,
\begin{equation}\label{eq:cs}
\xymatrix@-15pt{
B_0 &\subset & B_1 \\
\cup && \cup \\
A_0 & \ss & A_1
}
\end{equation}
along with a nondegenerate trace on $B_1$.  Given any inclusion of
algebras $A \subset B$ with a nondegenerate trace on $B$, let $E^B_A$
denote the unique trace preserving conditional expectations from $B$
to $A$.  The square (\ref{eq:cs}) is a
\emph{commuting} square if  $E^{B_1}_{A_1}(B_0)=A_0$.  The square is
\emph{symmetric} if $B_1$ is linearly spanned by $B_0A_1$.  There are
many equivalent conditions to these; see
\cite{ghj}, \cite{jones-sunder1} for details.  An inclusion of
finite-dimensional $C^\ast$--algebras $A \subset B$ is
\emph{connected} if its Bratteli diagram is connected (equivalently,
the centers of $B$ and $A$ have trivial intersection).  Assume
(\ref{eq:cs}) is a symmetric commuting square with connected
inclusions, and the trace is the unique Markov trace of the inclusion
$B_0 \ss B_1$.  We can construct a ladder of symmetric commuting
squares:
\begin{equation}
\xymatrix@-15pt{ 
B_0 &\ss &B_1 &\ss &B_2 &\ss &\dots &\ss &B_n &\ss &\dots \\
\cup  && \cup && \cup && \dots && \cup &&\dots \\
A_0 &\ss & A_1 &\ss &A_2 &\ss &\dots &\ss &A_n &\ss &\dots 
}
\end{equation}
by iterating the basic construction to the right, i.e. do the basic
construction to the right on the top row and adjoin the Jones
projection to the bottom row (see \cite{jones-sunder1}, \cite{ghj},
\cite{popa83:_orthog_neuman}).  Then we may complete the inclusion of
algebras $\cup_n A_n \ss
\cup_n B_n$ with respect to the unique trace on $\cup_n B_n$ to obtain
a hyperfinite ${\rm II}_1$ subfactor (\cite{jones-sunder1}).

Fix integers $k$ and $p$.  We will consider squares of the form
(\ref{eq:sq}).
Index the rows and columns of matrices in $M_p(\mathbb{C})\otimes
M_k(\mathbb{C})$ by the set $\{1,2, \ldots p \}\times \{1,2, \ldots ,
k\}$ in the natural way.  Following the notation in
\cite{krishnan-sunder1}, we denote elements of $\{1,2,\dots,p\}$ with
Greek letters and elements of $\{1,2,\dots,k\}$ with Roman letters.
For $F$ in $M_p(\mathbb{C})\otimes M_k(\mathbb{C})$, denote the entry
of $F$ in row $(\alpha, a)$ and column $(\beta, b)$ by $F_{\alpha
a}^{\beta b}$.

Let $U$ be a permutation matrix in $M_p(\mathbb{C})\otimes
M_k(\mathbb{C})$, i.e. its entries are
either $0$ or $1$.  $U$ is a permutation biunitary if 
the block transpose $\widetilde{U}$ of $U$, defined by 
\begin{equation}\label{eq:biunit}
\widetilde{U}_{\beta a}^{\alpha b} = U_{\alpha a}^{\beta b},
\end{equation}
is also a permutation matrix.  Equivalently, $U$ is a permutation
biunitary if 
\begin{equation}\label{eq:sq}
\xymatrix@-15pt{
U(1\otimes M_k(\mathbb{C}))U^\ast& \subset & M_p(\mathbb{C})\otimes
M_k(\mathbb{C}) \\
\cup && \cup \\
{\mathbb{C}} &\subset & M_p(\mathbb{C}) \otimes 1 
}
\end{equation}
is a commuting square
(\cite{jones-sunder1},\cite{y.90:_quant_symmet_differ_geomet_finit}).
We may then construct a subfactor as described above, which we denote
$R_U \ss R$.  In \cite{krishnan-sunder1}, Krishnan and Sunder list all
the nonequivalent subfactors of this type with $k=p=3$ and compute the
principal graphs of all the finite depth ones.

\section{The Intermediate Subfactor}
\label{sec:constr-interm-subf}
In Proposition~\ref{pro:intcomsq}, we show that if $U$ is a permutation
biunitary then the commuting square (\ref{eq:sq}) may be decomposed
into two adjacent symmetric commuting squares,

\begin{equation}\label{eq:sqs}
\xymatrix@-15pt{
U(1\otimes M_k(\mathbb{C}))U^\ast& \subset & M_p(\mathbb{C})\otimes
M_k(\mathbb{C}) \\
\cup && \cup \\
U(1\otimes \Delta_k )U^\ast & \subset & M_p(\mathbb{C}) \otimes \Delta_k \\
 \cup && \cup \\
{\mathbb{C}} & \subset & M_p(\mathbb{C}) \otimes 1 
}
\end{equation}
where $\Delta_k \cong \mathbb{C}^k$ is the diagonal subalgebra of
$M_k(\mathbb{C})$.  

An essential ingredient in our analysis of Krishnan-Sunder subfactors
is the following result of \cite{krishnan-sunder1}.

\begin{lem}[Krishnan-Sunder]\label{pro:ks}
If $U$ is a permutation biunitary in $M_p(\mathbb{C})\otimes
M_k(\mathbb{C})$, then there exist permutations
$\dgre_1,\dgre_2,\dots,\dgre_p$ in $S_k$ and permutations
$\dgge_1,\dgge_2,\dots,\dgge_k$ in $S_p$, such that
\[
U_{\alpha a}^{\beta b}= \delta_{\beta,\dgge_a(\alpha)}
\delta_{b,\dgre_\alpha(a)} 
\]
Consequently, there are permutations 
$\gre_1,\gre_2,\dots,\gre_p$ in $S_k$ and permutations
$\gge_1,\gge_2,\dots,\gge_k$ in $S_p$, such that
\[
(U^\ast)_{\alpha a}^{\beta b}= \delta_{\beta,\gge_a(\alpha)}
\delta_{b,\gre_\alpha(a)} 
\]
(In
Krishnan and Sunder's notation  \cite{krishnan-sunder1}, $\gre = \psi^{-1}$ and $\gge =
\phi^{-1}$).
\end{lem}

\begin{pro}\label{pro:intcomsq}
If $U$ is a permutation biunitary in $M_p(\mathbb{C})\otimes
M_k(\mathbb{C})$, then $U(1\otimes \Delta_k )U^\ast \subset
M_p(\mathbb{C}) \otimes \Delta_k$, and the two small squares in
(\ref{eq:sqs}) are symmetric commuting squares.  
\end{pro}
\begin{proof}
To simplify notation, in this proof we denote $\dgre_{\alpha}(a)$ and
$\dgge_a(\alpha)$ in Lemma~\ref{pro:ks} by $\alpha(a)$ and
$a(\alpha)$, respectively.  Let $\{ e_{\alpha,\beta}\}$ and
$\{f_{a,b}\}$ be the natural sets of matrix units for
$M_p(\mathbb{C})$ and $M_k(\mathbb{C})$, respectively.  Then $\{1
\otimes f_{a,a}\}_a$ is a basis of $1\otimes\Delta_k$, and
\begin{equation}\label{eq:intbasis}
U ( 1 \otimes f_{a,a}) U^\ast =
U \left( \sum_{\alpha=1}^p e_{\alpha,\alpha} \otimes f_{a,a}\right) U^\ast 
= \sum_{\alpha=1}^p e_{a(\alpha),a(\alpha)} \otimes
f_{\alpha(a),\alpha(a)} 
\end{equation}
which is contained in $M_p(\mathbb{C})\otimes \Delta_k$.
Hence, $U(1\otimes \Delta_k )U^\ast \subset M_p(\mathbb{C}) \otimes
\Delta_k$.  

Since (\ref{eq:sq}) is a symmetric commuting square, it suffices to
prove that the upper square is commuting and the lower square is symmetric.

The trace-preserving conditional expectation $E$ from $M_p(\mathbb{C})
\otimes M_k(\mathbb{C})$ to $M_p(\mathbb{C}) \otimes \Delta_k$, is given by
$$E(e_{\alpha,\beta} \otimes f_{a,b}) = \delta_{a,b}(e_{\alpha,\beta} \otimes f_{a,a})$$
So,
\begin{align*}
E\left(U ( 1 \otimes f_{a,b}) U^\ast \right) 
&= E\left(\sum_{\alpha=1}^p e_{a(\alpha),b(\alpha)} \otimes f_{\alpha(a),\alpha(b)}\right) \\
&= \delta_{a,b} \sum_{\alpha=1}^p(e_{a(\alpha),a(\alpha)} \otimes f_{\alpha(a),\alpha(a)}) \\
&= \delta_{a,b} U(1 \otimes f_{a,a})U^\ast
\end{align*}
Therefore, the upper square of (\ref{eq:sqs}) is commuting.  


It remains to prove that the lower square is symmetric.  For any
$b\in\{1,2,\dots, k\}$ and $\beta \in \{1,2,\dots, p\}$ there exists
$a, \alpha$ such that $U_{\alpha,a}^{\beta,b}=1$, because $U$ is a
permutation matrix.  By (\ref{eq:intbasis}) and the fact that $\dgge_a$
is a permutation,
\[
e_{a(\alpha),\gamma}\otimes f_{\alpha(a),\alpha(a)} \in
(U(1\otimes \Delta_k)U^\ast)(M_p(\mathbb{C})\otimes 1)
\]
for all $\gamma\in\{1,2,\dots,p\}$.  Since $\beta=a(\alpha)$ and
$b=\alpha(a)$ were chosen arbitrarily, $M_p(\mathbb{C})\otimes
\Delta_k$ is linearly spanned by $(U(1\otimes \Delta_k)U^\ast
)(M_p(\mathbb{C})\otimes 1)$.
\end{proof}

Note that the inclusion $U(1\otimes \Delta_k)U^\ast \ss
M_p(\mathbb{C})\otimes \Delta_k$ is not necessarily connected.  
Proposition~\ref{pro:intcomsq} immediately implies,
\begin{cor}\label{cor:dualsubf}
If $U$ is a permutation biunitary in $M_p(\mathbb{C})\otimes
M_k(\mathbb{C})$, then $R_U \ss R$ has an intermediate von Neumann
subalgebra $R_U \ss P_U \ss R$, i.e. the subalgebra constructed from
the lower symmetric commuting square in (\ref{eq:sqs}).  In
particular, if $R_U \ss R$ is irreducible then $R_U \ss P_U \ss R$ is
an intermediate subfactor.
\end{cor}

Recall $\widetilde{U}$ in (\ref{eq:biunit}).  We will show in
Proposition~\ref{pro:dual} that $R_{\widetilde{U}} \ss
P_{\widetilde{U}}$ and $P_U \ss R$ are dual (by symmetry so are $R_{U}
\ss P_{U}$ and $P_{\widetilde{U}} \ss R$) and can be constructed from
a biunitary permutation matrix $\myA\in\Delta_p\otimes M_k(\mathbb{C})$

Let the permutations $\dgre_\alpha$, $\dgge_a$, $\gre_\alpha$, and
$\gge_a$ be as in Lemma~\ref{pro:ks}.  Define permutation matrices $\myA\in
\Delta_p\otimes M_k(\mathbb{C})$ and $\myB\in M_p(\mathbb{C})\otimes
\Delta_k$ by:
\[
\myA_{\alpha a}^{\beta b} = \delta_{\beta, \alpha}\delta_{b, \dgre_{\alpha}(a)}
\qquad
\myB_{\alpha a}^{\beta b} = \delta_{\beta, \dgge_a(\alpha)}\delta_{b, a}
\]
and permutation matrices
$\myC\in \Delta_p\otimes M_k(\mathbb{C})$ and $\myD\in
M_p(\mathbb{C})\otimes \Delta_k$ by:
\[
\myC_{\alpha a}^{\beta b} = \delta_{\beta, \alpha}\delta_{b, \gre_{\alpha}(a)}
\qquad
\myD_{\alpha a}^{\beta b} = \delta_{\beta, \gge_a(\alpha)}\delta_{b, a}.
\]
It is easy to check that
\begin{equation}\label{eq:uabc}
U=\myD^\ast \myA = \myC^\ast \myB
\end{equation}

and

\begin{equation}\label{eq:duabc}
\widetilde{U}= \myD \myC^\ast = \myA \myB^\ast
\end{equation}

The following lemma shows that we may replace $U$ in the commuting
squares that engender $R_U \ss P_U$ and $P_U \ss R$, by the simpler
unitaries $\myC^\ast$ and $\myA$.
\begin{lem}\label{lem:intsq}
Let $\myA,\myC \in \Delta_p\otimes M_k(\mathbb{C})$ and $\myB,\myD \in
M_p(\mathbb{C})\otimes \Delta_k$ be permutation matrices satisfying
(\ref{eq:uabc}) and (\ref{eq:duabc}).  Then

\begin{equation}\label{eq:lowsq2}
\xymatrix@-15pt{ U(1\otimes \Delta_k )U^\ast &
\subset & M_p(\mathbb{C}) \otimes \Delta_k \\ 
\cup && \cup \\
{\mathbb{C}} &\subset & M_p(\mathbb{C}) \otimes 1 }
\qquad = \qquad
\xymatrix@-15pt{ \myC^\ast (1\otimes \Delta_k )\myC &
\subset & M_p(\mathbb{C}) \otimes \Delta_k \\ 
\cup && \cup \\
{\mathbb{C}} &\subset & M_p(\mathbb{C}) \otimes 1 }
\end{equation}

and 

\begin{equation}\label{eq:upsq2}
\xymatrix@-15pt{ U(1 \otimes M_k({\mathbb{C}} ))U^\ast &
\subset & M_p(\mathbb{C}) \otimes M_k(\mathbb{C}) \\ 
\cup && \cup \\
U(1\otimes \Delta_k )U^\ast &
\subset & M_p(\mathbb{C}) \otimes \Delta_k \\ 
}
\qquad\cong\qquad
\xymatrix@-15pt{ 
\myA (1\otimes M_k(\mathbb{C}) ) \myA^\ast & \subset & M_p(\mathbb{C})
\otimes M_k(\mathbb{C}) \\  
\cup && \cup \\
\myA (1\otimes \Delta_k )\myA^\ast & \subset & M_p(\mathbb{C}) \otimes \Delta_k  }
\end{equation}

\end{lem}
\begin{proof}
Note that conjugation by any permutation matrix in $M_p(\mathbb{C})\otimes
\Delta_k$ (such as $\myB$) stabilizes $1\otimes \Delta_k$.
So substituting $\myC^\ast\myB$ for $U$ yields (\ref{eq:lowsq2}).
Similarly, substitute $U=\myD^\ast\myA$ in left-hand side of (\ref{eq:upsq2}),
and then conjugate the entire square by $\myD$ to obtain the right-hand side.
\end{proof}

\begin{pro}\label{pro:dual}
$P_U \ss R$ and $R_{\widetilde{U}} \ss
P_{\widetilde{U}}$ are dual inclusions.
\end{pro}
\begin{proof}
Clearly $\widetilde{\widetilde{U}} = U$.
So by symmetry, Lemma~\ref{lem:intsq} implies
\begin{equation}\label{eq:dlowsq}
\xymatrix@-15pt{ \widetilde{U}(1\otimes \Delta_k )\widetilde{U}^\ast &
\subset & M_p(\mathbb{C}) \otimes \Delta_k \\ 
\cup && \cup \\
{\mathbb{C}} &\subset & M_p(\mathbb{C}) \otimes 1 }
\qquad \cong \qquad
\xymatrix@-15pt{ \myA (1\otimes \Delta_k )\myA^\ast &
\subset & M_p(\mathbb{C}) \otimes \Delta_k \\ 
\cup && \cup \\
{\mathbb{C}} &\subset & M_p(\mathbb{C}) \otimes 1 }
\end{equation}
Let $e\in M_k(\mathbb{C})$ be the usual Jones projection of the Jones
extension $\mathbb{C} \ss \Delta_k \ss M_k(\mathbb{C})$.  Conjugation
by $\myA$ fixes $1\otimes e$, hence (\ref{eq:upsq2}) is the upward
basic construction of (\ref{eq:dlowsq}).  Since both squares are
symmetric, doing the basic contruction to the right yields dual
inclusions.
\end{proof}

\section{The Subgroup Subfactor}\label{ss:subgroup} 
Fix a permutation biunitary $U$.
Let $\gre_{\alpha}$, $1\leq \alpha \leq p$, and
$\gge_{a}$, $1\leq a \leq k$, be the permutations defined before
Lemma~\ref{lem:intsq}.  

Let $\grp$ be the subgroup of $S_k$ generated by elements of the form 
$\gre_{\alpha} \gre_\beta^{-1}$, i.e., 
\[
\grp =\grp_U= \langle \gre_{\alpha} \gre_\beta^{-1} :
\alpha,\beta\in\{1,2,\dots, p\}\rangle 
\]
For each $a\in \{1,2,\dots,k\}$ let $H_a\ss \grp$ be the subgroup that
fixes $a$.
If $\grp$ acts transitively on $\{1,2,\dots,k\}$, then the subgroups
$H_a$ are all conjugate.  In this case set $H= H_1$.  In general, let
$\Omega$ be the set of orbits in $\{1,2,\dots, k\}$, and for each 
$r\in \Omega$ set $H_r = H_a$ for an arbitrary representative $a$ in
$r$.

\begin{rem}\label{rem:gammagroup}
Krishnan and Sunder use the group $\langle \gre_{\alpha}^{-1}
\gre_\beta: \alpha,\beta\in\{1,2,\dots, p\}\rangle$ instead of $\grp$
in \cite{krishnan-sunder1}.  However, the two groups, as well as their
fixed point subgroups, are equivalent via conjugation by $\gre_\gamma$
for any $\gamma\in\{1,2,\dots, p\}$.  We depart from
\cite{krishnan-sunder1} for notational convenience in the proof of the
following theorem.
\end{rem}

\begin{thm}\label{thm:ks-subgr-iso}
Let $\grp$ act on $R$ by outer automorphisms.  There is a canonical
isomorphism of $Z(P_U)$ with $\mathbb{C}^{\Omega}$.  If $q_r$ is the minimal
projection in $Z(P_U)$ corresponding to $r\in\Omega$, then $q_r
R_U \ss q_r P_U $ is isomorphic to $R^{\grp} \ss R^{H_r}$.  In
particular, if $\grp$ acts transitively on $\{1,2,\dots,k\}$ then $R_U
\ss P_U$ is a subfactor, and $R_U \ss P_U$ is isomorphic to $R^{\grp} \ss
R^H$.
\end{thm}
\begin{proof}
By construction of $\grp$, the cosets $\gre_{\alpha}^{-1}\grp$ 
are the same for all $\alpha\in\{1,2,\dots,p\}$.  Denote this coset by
$\dgrp$.  
Given a set $S$, let $\Delta_S$ denote the algebra of functions
$S\rightarrow \mathbb{C}$ with pointwise multiplication and the
trace $f \mapsto \frac{1}{\vert S \vert}\sum_{s\in S} f(s)$.  Denote
the characteristic function of $s\in S$ by $x_s$.
Define an inclusion map
$i\colon \Delta_\grp \to M_p(\mathbb{C})\otimes \Delta_{\dgrp}$ by
$i(x_g) = \sum_{\alpha} e_{\alpha, \alpha} \otimes x_{\gre_{\alpha}^{-1} g }$.  
It is straightforward to check that
\begin{equation}\label{eq:grpsq}
\xymatrix@-15pt{ 
\Delta_\grp  &\subset_i & M_p(\mathbb{C})\otimes \Delta_{\dgrp} \\  
\cup && \cup \\
\mathbb{C} & \subset & M_p(\mathbb{C}) }
\end{equation}
is a symmetric commuting square with connected inclusions, if we take
the trace on $M_p(\mathbb{C})\otimes \Delta_{\dgrp}$ to be the product
trace.

Then we can construct a hyperfinite $\rm{II}_1$ subfactor $\salg \ss
\alg$ by iterating the basic construction to the right in the usual way.
\begin{equation}\label{eq:bc}
\xymatrix@-15pt{ 
\Delta_\grp  &\subset^\grph & M_p(\mathbb{C})\otimes \Delta_{\dgrp} 
&\ss^{\grph^t} &\alg_2 &\ss^\grph &\alg_3 &\ss^{\grph^t} &\dots   &\alg\\  
\cup && \cup && \cup && \cup &&& \cup\\
\mathbb{C} & \subset & M_p(\mathbb{C}) 
&\ss &\salg_2 &\ss &\salg_3 &\ss &\dots   &\salg \\  }
\end{equation}
Note that the Bratteli diagram $\grph$ as marked in (\ref{eq:bc}) is
the bipartite graph with even vertices labeled by $\grp$, odd vertices
labeled by $\dgrp$, and an edge for each pair $(g, \alpha)\in
\grp\times\{1,2,\dots,p\}$ going from $g$ to $\gre_{\alpha}^{-1} g$.  We
denote the reflection of $\grph$ by $\grph^t$.  For each $n$, $\salg_n
\ss \alg_n$ is isomorphic to $M_{p^n}(\mathbb{C})\otimes 1 \ss
M_{p^n}(\mathbb{C})\otimes \Delta_{\tilde{\grp}} $, where
$\tilde{\grp}$ is $\grp$ or $\dgrp$ according to whether $n$ is even
or odd.  We claim that $\salg \ss \alg$ is irreducible.  By Ocneanu
compactness, $\salg' \cap \alg = M_p(\mathbb{C})' \cap \Delta_{\grp} =
\Delta_{\grp} \cap (1 \otimes \Delta_{\dgrp})$.
Suppose $\sum_{g\in\grp} k_g x_g \in 1 \otimes \Delta_{\dgrp}$, where $k_g
\in \mathbb{C}$.  Then 
\[
\sum_{g\in\grp} k_g x_g 
= \sum_{g'\in\dgrp} \sum_{\gre_{\alpha}^{-1}g = g'} k_g (f_{\alpha,
\alpha} \otimes x_{g'})
\]
Since $\sum k_g x_g \in 1 \otimes \Delta_{\dgrp}$, we must have that
$k_g$ is constant over all $g \in \{\gre_\alpha g'\}_\alpha$.  Since
this is holds for all $g'$, it follows that $k_g$ is constant over all
$g \in \grp 1 = \grp$.  Therefore, $\sum k_g x_g \in \mathbb{C}$.
This proves the claim.  

For each $g\in \grp$, let $\mu_g$ be the automorphism of $\grph$ that
maps each vertex $g'\in \grp\cup\dgrp$ to $g'g^{-1}$ and each edge
$(g',\alpha)$ to the edge $(g'g^{-1}, \alpha)$.  The morphism is well
defined since the endpoints of the edge $(g', \alpha)$ are mapped to
the endpoints of $(g'g^{-1}, \alpha)$.  Moreover, $\mu_g$ obviously
preserves the trace weights of $\grph$.  Clearly $g \mapsto \mu_g$ is
a group action of $G$ on $\grph$.  Now extend $\mu$ to the chain of
Bratteli diagrams of the top row of inclusions in (\ref{eq:bc}).  For
each $n$, $\mu$ implements an action $\mu^n$ of $\grp$ on $\alg_n$ by
trace preserving automorphisms.  The family of actions $\{\mu^n:
\grp\rightarrow \aut(\alg_n)\}_n$ are consistent,
i.e. $\mu^n\vert_{\alg_{n-1}} = \mu^{n-1}$, and thus extend to an
action of $\grp$ on $\alg$.  We denote this action again by $\mu$.
Note that the action of $\grp$ on $\alg$ is outer since $\salg \ss
\alg$ is irreducible.  Therefore $\salg \ss \alg$ is isomorphic to
$R^\grp \ss R$, as in the statement of the theorem.

Let $E^{\grp}$ be the group averaging maps from $\alg$
onto the fixed point algebra $\alg^\grp$.
The action of $\mu_g$ on $\alg_n = M_{p^n}(\mathbb{C})\otimes
\Delta_{\tilde{\grp}}$ is given by $\mu_g(F \otimes x_g') = F\otimes
x_{g'g^{-1}}$, hence $E^{\grp}(\alg_n)=\salg_n$ for each $n$.
Therefore, $\alg^{\grp}=\salg$.

We first assume that $\grp$ acts transitively on $\{1,2,\dots, k\}$.
Define an inclusion $\Delta_{\grp / H} \to \Delta_{\grp}$ by
$x_{gH} \mapsto \sum_{g'\in gH} x_{g'}$.  Similarly define
$\Delta_{\dgrp / H} \to \Delta_{\dgrp}$.  For each $n$, let 
$\ialg_n = \alg_n^H$.  Clearly $\ialg_n = \salg_n \otimes
\Delta_{\tilde{\grp}/H}$ where $\tilde{\grp}$ is $\grp$ or $\dgrp$
according to whether $n$ is even or odd. Thus we have an intermediate
chain of algebras
\begin{equation}\label{eq:ibc}
\xymatrix@-15pt{ 
\Delta_\grp  &\subset^\grph & M_p(\mathbb{C})\otimes \Delta_{\dgrp} 
&\ss^{\grph^t} &\alg_2 &\ss^\grph &\alg_3 &\ss^{\grph^t} &\dots   &\alg\\  
\cup && \cup && \cup && \cup &&& \cup\\
\Delta_{\grp/H}  &\subset & M_p(\mathbb{C})\otimes \Delta_{\dgrp/H} 
&\ss &\ialg_2 &\ss &\ialg_3 &\ss &\dots
&\ialg\\  
\cup && \cup && \cup && \cup &&& \cup \\
\mathbb{C} & \subset & M_p(\mathbb{C}) 
&\ss &\salg_2 &\ss &\salg_3 &\ss &\dots   &\salg \\  }
\end{equation}
where $\ialg = \alg^H$.  Since the group averaging map $E^H$ from $A$
onto $A^H$ is the conditional expectation from $\alg_n$ onto $\ialg_n$
for each $n$, it follows that the upper-left-most square of
(\ref{eq:ibc}) is commuting.  It is straightforward to verify that the
lower-left-most square of (\ref{eq:ibc}) is symmetric, hence both of
the left-most squares are symmetric commuting squares.  For $n \geq
2$, $\ialg_n$ contains the Jones projection of the inclusion
$\alg_{n-2} \ss \alg_{n-1}$, hence the chain $(\ialg_n)_n$ contains
the Jones tower of $\ialg_0 \ss \ialg_1$.  Moreover, the Bratteli
diagram of $\ialg_{n-1} \ss \ialg_{n}$ is the transpose of
$\ialg_{n-2} \ss \ialg_{n-1}$ for $n \geq 2$, hence by dimension
considerations, the chain $(\ialg_n)_n$ is no more than the Jones
tower.  Therefore, both the upper and the lower ladders are the ones
obtained by iterating the basic construction in the usual way from the
left-most square.

Now consider the lower left square ($\ast$) of (\ref{eq:ibc}).
We claim that ($\ast$) is isomorphic to (\ref{eq:lowsq2}) via the
identification of $\grp/ H$ and $\dgrp/ H$ with $\{1,2,\dots, k\}$, by
$gH \mapsto g(1)$.  Let $(\hat{a})_{a=1}^k$ be the minimal projections in
$\Delta_k$, and define an isomorphism $M_p(\mathbb{C})\otimes
\Delta_{\dgrp/H} \rightarrow M_p(\mathbb{C})\otimes \Delta_{k}$ by
$F\otimes x_{gH} \mapsto F\otimes \hat{g(1)}$.  Fix $a$ and choose $f\in
\grp$ such that $f(1)=a$.  Then $\myC^\ast (1\otimes \hat{a}) \myC =
\sum_{\alpha} e_{\alpha,\alpha} \otimes (\gre_{\alpha}^{-1}(a))\hat{}$
is the image of $\sum_{\alpha} e_{\alpha,\alpha} \otimes
x_{\gre_{\alpha}^{-1}fH} = 1\otimes x_{fH}$.  Therefore, ($\ast$) is
isomorphic to (\ref{eq:lowsq2}), and $R_U \ss P_U$ is isomorphic to
$\alg^{\grp} \ss \alg^{H}$.  This proves $R_U \ss P_U$ is isomorphic
to $R^{\grp} \ss R^{H}$, as in the statement of the theorem.

If $\grp$ does not act transitively, then $\myC^\ast (1\otimes \Delta_k)
\myC \ss M_p(\mathbb{C})\otimes \Delta_k$ is not a connected
inclusion;  
its connected components correspond to the orbits
of $\grp$ in $\{1,2,\dots, k\}$.  
Given an orbit $r\in \Omega$, let $q_r = \myC^\ast \sum_{a\in r} \hat{a} \myC$.  
Clearly $q_r$ is central in $P_U$ and $q_r R_U \ss q_r P_U$ can be
obtained by iterating the basic construction on 
\begin{equation*}
\xymatrix@-15pt{ q_r\myC^\ast (1\otimes \Delta_k )\myC &
\subset & q_r (M_p(\mathbb{C}) \otimes \Delta_k) \\ 
\cup && \cup \\
q_r{\mathbb{C}} &\subset & q_r(M_p(\mathbb{C}) \otimes 1) }
\end{equation*}
By an identical argument as the one above (using the group $H_r$
instead of $H$), $q_r R_U \ss q_r P_U$ is isomorphic to $\alg^{\grp} \ss
\alg^{H_r}$.  This proves $q_r R_U \ss q_r P_U$ is isomorphic to
$R^{\grp} \ss R^{H_r}$, as in the statement of the theorem.  Then,
$q_r R_U \ss q_r P_U$ is a subfactor for each $r$, which implies that
$Z(P_U) = \bigoplus_{r\in\Omega} \mathbb{C}q_r$.

%
%
\end{proof}

\begin{cor}\label{cor:subgroup-ks}
Let $H\ss \grp$ be any inclusion of finite groups, and let $\grp$ act
on the hyperfinite $\rm{II}_1$ factor $R$ by an outer action.  Let
$k=\vert \grp/H \vert$.  Suppose the action of $\grp$ on $\grp/H$ can
be generated by $p'$ elements of $\grp$.  Then there exists a
permutation biunitary $U\in \Delta_{p'+1}\otimes M_k(\mathbb{C})$ such
that $R_U
\ss P_U$ (as defined in Corollary~\ref{cor:dualsubf}) is isomorphic to
$R^{\grp} \ss R^H$.
\end{cor}
\begin{proof}
Let $U_1, U_2, \dots U_{p'}$ be $k\times k$ permutation matrices that
generate the action of $\grp$ on $\grp/H$.  Set $U_0$ to be the
identity matrix.  Let $\{e_{\alpha}\}_{0\leq \alpha \leq p'}$ be a
basis of $\Delta_{p'+1}$, and set $U=\sum_{0}^{p'} e_{\alpha}\otimes
U_\alpha$.  Obviously $\grp_U = \grp$ and the fixed point subgroup of
$\grp$'s action on $\{1,2,\dots, k\}$ is $H$.  By
Theorem~\ref{thm:ks-subgr-iso}, this $U$ does the job.
\end{proof}

\section{The Bisch Projection}
\label{sec:comp-stand-invar}
We now show that the Bisch projection defined in
\cite{bhattacharyya98:_krish_sunder} (see also
\cite{binazeph:_inter_subfac}) corresponds to the intermediate
sub-von Neumann algebra $R_U \ss P_U$.  

The upward basic construction of $M_p(\mathbb{C})\otimes 1 \ss
M_p(\mathbb{C})\otimes M_k(\mathbb{C})$ in (\ref{eq:sq}) is
$M_p(\mathbb{C})\otimes M_{k^2}(\mathbb{C})$.  The matrix rows and
columns of its subalgebra $1 \otimes M_{k^2}(\mathbb{C})$ are indexed
naturally by the set $\{1,2,\dots,k\}\times\{1,2,\dots,k\}$
(\cite{y.90:_quant_symmet_differ_geomet_finit},\cite{jones-sunder1}).
Given $x \in 1 \otimes M_{k^2}(\mathbb{C})$, denote by $x_{ab}^{cd}$
the entry of $x$ in row $\{a,b\}$ and column $\{c,d\}$.  The first
relative commutant of $R_U \ss R$ is the subalgebra of $1 \otimes
M_{k^2}(\mathbb{C})$ satisfying the Ocneanu compactness condition
\cite{y.90:_quant_symmet_differ_geomet_finit}.  

Define $p\in 1\otimes M_{k^2}(\mathbb{C})$ by
\begin{equation}\label{eq:pdef}
p_{ab}^{cd} = 
\begin{cases}
1, &\text{if $a=b=c=d$}; \\
0, &\text{otherwise}
\end{cases}
\end{equation}

\begin{pro}\label{pro:bischproj0}
The projection $p$ defined above is contained in the first relative
commutant of $R_U \ss R$; and $p$ is a Bisch projection, that is, pnd$p$
implements the conditional expectation from $R$ to $\{p\}' \cap R$
with respect to the trace.
\end{pro}
\begin{proof}
This is proved in somewhat different notation in Lemma 6.4.1 of
\cite{bhattacharyya98:_krish_sunder}.  
For the convenience of the reader we give a proof here.  We first show that
$p$ is in the first relative commutant using Jones' diagrammatic
formulation of the higher relative commutants of $R_U' \cap R$
\cite{jones-sunder1}.

We claim that for $a, b \in \{1,2,\dots,k\}$
and $\alpha, \beta \in \{1,2,\dots,p\}$:
\begin{equation}\label{eq:pdiagram}
\vcenter{\xymatrix@R-5pt@C-12pt{
& b & b' \ar'[d][dd] &    \\
\beta \ar[rrr] & & & \alpha \\
& a \ar'[u][uu] & a' & \\
}}
\quad
= \quad
\begin{cases}
\delta_{\alpha, \beta}\cdot\delta_{b, b'} \cdot \delta_{a, \nu_{\beta}(b)}, &\text{if $a=a'$} \\
\delta_{\alpha, \beta}\cdot\delta_{a, a'} \cdot \delta_{a, \nu_{\beta}(b)}, &\text{if $b=b'$}
\end{cases}
\end{equation}
The claim is obvious from Section~5 of
\cite{bhattacharyya:_group_action_graph}, but here is a direct proof.
The left-hand side of (\ref{eq:pdiagram}) is by definition
$\sum_{\gamma \in \{1,2,\dots,p\}} U_{\gamma a}^{\beta b}
\overline{U_{\gamma a}^{\alpha b}}$.  Note that in our case entries of
$U$ are either 0 or 1, so $\overline{U} = U$.  We have:
\begin{align*}
\sum_{\gamma \in \{1,2,\dots,p\}} U_{\gamma a}^{\beta b} \overline{U_{\gamma a'}^{\alpha b'}}
&= \sum_{\gamma} \delta_{\beta, \rho_a(\gamma)} \delta_{b,\lambda_{\gamma}(a)} 
                 \delta_{\alpha, \rho_{a'}(\gamma)} \delta_{b', \lambda_{\gamma}(a')} \\
&= \delta_{\alpha, \beta}\delta_{b, b'}
\delta_{b,\lambda_{\rho_{a}^{-1}(\alpha)}(a)} &\text{if $a=a'$} \\
&= \delta_{\alpha, \beta}\delta_{b, b'}\delta_{a, \nu_{\alpha}(b)}
&\text{(Lemma~5 of \cite{krishnan-sunder1})} 
\end{align*}
and
\begin{align*}
\sum_{\gamma \in \{1,2,\dots,p\}} U_{\gamma a}^{\beta b} \overline{U_{\gamma a'}^{\alpha b'}}
&= \sum_{\gamma} (U^\ast)_{\beta b}^{\gamma a}(U^\ast)_{\alpha b'}^{\gamma a'} \\
&= \sum_{\gamma} \delta_{\gamma, \theta_b(\beta)} \delta_{a,\nu_{\beta}(b)} 
                 \delta_{\gamma, \theta_{b'}(\alpha)} \delta_{a', \nu_{\alpha}(b')} \\
&= \delta_{\alpha, \beta}\delta_{a, a'} \delta_{a, \nu_{\beta}(b)} &\text{if $b=b'$}
\end{align*}
This proves the claim.  
Using (\ref{eq:pdiagram}), it is easy to verify that $p$ satisfies the
diagrammatic condition for $p$ to be in $R_U' \cap R_1$ (see Theorem
6.1.4 and the preceding discussion in \cite{jones-sunder1}).  

Let $q \in M_{k^{2}}(\mathbb{C})$ be the projection identified with
$p$, that is, $p = 1 \otimes q \in M_p(\mathbb{C}) \otimes
M_{k^{2}}(\mathbb{C})$.  
Define $p_n = 1 \otimes 1_{M_{k^n}(\mathbb{C})} \otimes q \in
M_p(\mathbb{C}) \otimes M_{k^{n+2}}(\mathbb{C})$, for $n = 0, 1, 2,
\dots$.  The same argument as above shows that $p_n$ is contained in
$(n-1)$st relative commutant of $R_U \ss R$.
It is easy to verify that the sequence of projections $(p_n)$ along
with the sequence of Jones projections $(e_n)$ of $R_U \ss R$
satisfies the Bisch-Jones relations of the Fuss-Catalan algebras
\cite{bisch-jones1} with $\alpha = \beta = k$.  It follows by
\cite{MR95c:46105} that $p$ implements the conditional expectation
from $R$ to $\{p\}' \cap R$ with respect to the trace
\end{proof}

\begin{pro}\label{pro:bischproj}
Let $p$ be the Bisch projection defined above and let $E_p$ be the
conditional expection from $R$ to $\{p\}' \cap R$ implemented by $p$
(by Proposition~\ref{pro:bischproj0}).  Then $P_U = E_p(R)$.
\end{pro}
\begin{proof}
Let $e$ be the Jones projection of the extension $R_U \ss R \ss R_1$.
Let $(A_m)$, $(C_m)$, and $(B_m)$, $m = 0,1,2,\dots$, be the
chains of algebras obtained by iterating the basic construction on
(\ref{eq:sqs}) to the right, and let $D_m = \langle B_m, e \rangle$:
\begin{equation}
\xymatrix@-15pt{ 
\calg_0  &\subset & \calg_1
&\ss &\calg_2 &\ss&\calg_3 &\ss&\dots   & R_1\\  
\cup && \cup && \cup && \cup &&& \cup\\
U(1\otimes M_k(\mathbb{C}))U^\ast  &\subset & M_p(\mathbb{C})\otimes M_k(\mathbb{C})
&\ss &\alg_2 &\ss&\alg_3 &\ss&\dots   & R\\  
\cup && \cup && \cup && \cup &&& \cup\\
U(1\otimes \Delta_k)U^\ast  &\subset & M_p(\mathbb{C})\otimes \Delta_k
&\ss &\palg_2 &\ss &\palg_3 &\ss &\dots
& P_U\\  
\cup && \cup && \cup && \cup &&& \cup \\
\mathbb{C} & \subset & M_p(\mathbb{C}) 
&\ss &\salg_2 &\ss &\salg_3 &\ss &\dots   & R_U \\  }
\end{equation}
The type of construction above is well-known
(\cite{y.90:_quant_symmet_differ_geomet_finit}, also see Ocneanu
compactness in \cite{jones-sunder1}), and in particular has the
property that the chain $\salg_m \ss \palg_m \ss \alg_m \ss \calg_m$ 
has the same Bratteli diagram for all odd $m$.  In other words, $\salg_m
\ss \palg_m \ss \alg_m \ss \calg_m$ is isomorphic to $\salg_m
\ss \salg_m \otimes \Delta_k \ss \salg_m \otimes M_{k}(\mathbb{C})
\ss \salg_m \otimes M_{k^2}(\mathbb{C})$ for all odd $m$.

By definition, $p$ is a projection in $1 \otimes M_{k^2}(\mathbb{C}) \ss
\calg_1 = M_p(\mathbb{C})\otimes M_{k^2}(\mathbb{C})$.  
Let $q \in M_{k^2}(\mathbb{C})$ be the projection identified with $p$,
that is, $p = 1 \otimes q \in M_p(\mathbb{C})\otimes
M_{k^2}(\mathbb{C})$.  Using (\ref{eq:pdef}), it is easy to check:
\begin{equation}\label{eq:qprop} 
q \in (\Delta_k \otimes 1)' \cap M_{k^2}(\mathbb{C})
\end{equation}
\begin{equation}\label{eq:qprop1} 
q (M_{k}(\mathbb{C}) \otimes 1) q = (\Delta_k \otimes 1)q.
\end{equation}
It follows that $p \in \palg_1' \cap \calg_1$ and $p \alg_1 p =
\palg_1 p$, hence $\palg_1 = E_p(\alg_1)$.  

By Proposition~\ref{pro:bischproj0}, $p$ is contained in $
R_U' \cap R_1$, hence $p$
is \emph{flat} \cite{y.90:_quant_symmet_differ_geomet_finit}.  By
flatness, $p$, as an element of $\salg_m \otimes M_{k^2}(\mathbb{C}) \cong
\calg_m$, for $m$ odd, is identified with $1_{\salg_m} \otimes q \in \salg_m
\otimes M_{k^2}(\mathbb{C})$.  Recalling that $\salg_m
\ss \palg_m \ss \alg_m \ss \calg_m$ is isomorphic to $\salg_m
\ss \salg_m \otimes \Delta_k \ss \salg_m \otimes M_{k}(\mathbb{C})
\ss \salg_m \otimes M_{k^2}(\mathbb{C})$, we have by (\ref{eq:qprop})
and (\ref{eq:qprop1}) that $p \in \palg_m' \cap \calg_m$ and $p \alg_m
p = \palg_m p$ for all odd $m$.  Hence, $\palg_m = E_p(\alg_m)$ for
all odd $m$.  Then by weak continuity of $E_p$, $P_U = E_p(R)$.  
\end{proof}

We restate Proposition~\ref{pro:bischproj} and
Proposition~\ref{thm:ks-subgr-iso} as follows:
\begin{cor}
The intermediate sub-von Neumann algebra of $R_U \ss R$ corresponding
to the Bisch projection $p$ is $R_U\ss P_U$, as defined in
Corollary~\ref{cor:dualsubf}.
\end{cor}

\end{document}